\documentclass[a4paper,11pt,twoside]{article}
\usepackage{pst-fill,pst-grad}
\usepackage{textcomp}
\usepackage[english]{babel}
\usepackage[utf8]{inputenc} 
\usepackage{graphicx}
\usepackage{amsmath}
\usepackage{float}
\usepackage{fancyhdr}
\usepackage[matrix,arrow,curve]{xy}
\usepackage{pstricks} 
\usepackage{amsmath,amsfonts,verbatim,afterpage,theorem,euscript,mathrsfs,amssymb}
\usepackage{amsfonts}
\usepackage{amssymb}
\usepackage{array}
\usepackage{dsfont}
\usepackage{hyperref}
\usepackage{authblk}

\newtheorem{Proposition}{Proposition}[section]

\newtheorem{Theoreme}{Theorem}
\newtheorem{Corollaire}{Corollary}

\def \P{\mathbb{P}}
\def \U{\vec{U}}
\def \W{\vec{W}}

\def \Rt{\mathbb{R}^{3}}

\def \finpv{\hfill $\blacksquare$  \\ \newline }
\def \pv{{\bf{Proof.}}~}

\def \ds{\displaystyle}
\title{ \bf A remark on the Liouville problem for stationary Navier-Stokes equations in Lorentz and Morrey spaces}
\author[1]{\small OSCAR JARR\'IN\footnote{or.jarrin@uta.edu.ec}}
\affil[1]{\scriptsize Dirección de investigación y desarrollo (DIDE).
Universidad Técnica de Ambato, campus Huachi, 
Avenida de los Chasquis y rio Payamino, 180207, Ambato, Ecuador.}
\begin{document} 
\maketitle
\begin{scriptsize}
\abstract{The Liouville problem for the stationary Navier-Stokes equations on the whole space is a challenging open problem who has know several recent contributions.  We prove here some  Liouville type theorems for these equations   provided the velocity field belongs to some Lorentz spaces  and then in the  more general setting of Morrey spaces. Our theorems correspond to a improvement of some recent results on this problem and contain some well-known results as a particular case.}\\[5mm]
\textbf{Keywords: Navier--Stokes equations; stationary system; Liouville theorem; Lorentz spaces; Morrey spaces}
\end{scriptsize}
\section{Introduction} 
In this article we review some recent results on the  Liouville problem for the stationary and incompressible Navier-Stokes equations in the whole space $\Rt$: 
\begin{equation}\label{N-S-f-nulle-intro}
-\Delta \U+(\U \cdot \vec{\nabla})\U +\vec{\nabla}P = 0, \qquad div(\U)=0, \\
\end{equation}
where $\U:\Rt\longrightarrow \Rt$ is the velocity and $P:\Rt\longrightarrow \mathbb{R}$ is the pressure. Recall that a weak solution of these equations is a couple $(\U, P) \in L^{2}_{loc}(\Rt)\times \mathcal{D}'(\Rt)$. Moreover,  since the pressure $P$ is always related to the velocity $\U$ by the identity $P=\frac{1}{-\Delta}\left( div\big((\U \cdot \vec{\nabla})\U\big ) \right)$ then we can concentrate our study in the variable $\U$. \\

The classical Liouville problem for the stationary Navier-Stokes equations states that the unique solution of equations (\ref{N-S-f-nulle-intro}) which verifies 
\begin{equation*}
\int_{\Rt} \vert \vec{\nabla} \otimes \U (x) \vert^2 dx <+\infty,
\end{equation*}
and 
\begin{equation*}
\vert \U(x)\vert \longrightarrow 0, \quad \text{as}\quad \vert x \vert \longrightarrow +\infty,
\end{equation*}
is the trivial solution $\U=0$, see the book \cite{Galdi},  the PhD thesis \cite{Jarrin} and the articles \cite{ChaeYoneda,ChaeWolf,ChaeWeng,Ser2015,Ser2016} for more references. Even though an answer to this question is not yet available,  great efforts have been invested to understand this open problem. More precisely, the main idea is to give some \emph{a priori} conditions on the decaying of solution $\U$  which allow us to prove that  $\int_{\Rt} \vert \vec{\nabla} \otimes \U (x) \vert^2 dx<+\infty$, and  with this information at hand, and sometimes with supplementary hypothesis on the solution $\U$,  we look for the identity $\U=0$. \\

In this setting, one of the first results is due to G. Galdi, see Theorem $X.9.5$ (page 729 ) of the book \cite{Galdi}, where it is proven that  if  $\U \in L^{\frac{9}{2}}(\Rt)$ then we have $\int_{\Rt} \vert \vec{\nabla} \otimes \U \vert^{2}dx \leq c \Vert \U \Vert^{3}_{L^{\frac{9}{2}}}$, and moreover, it is proven  the following local estimate for all $R>0$ and $c>0$ a constant independent of $R$:
$$ \int_{B_{R/2}} \vert \vec{\nabla} \otimes \U \vert^2 dx \leq c \Vert \U \Vert^{3}_{L^{\frac{9}{2}} (\mathcal{C}(R/2,R))},$$ where  $B_R =\{x\in \Rt: \vert x \vert < R \}$ and  $\mathcal{C}(R/2, R)=\{x\in \Rt: \frac{R}{2} \leq \vert x \vert \leq R\}$, which   yields the identity $\U=0$ provided  that $\U \in L^{\frac{9}{2}}(\Rt)$. \\ 

Galdi's result was thereafter extend to the  Lorentz space $L^{\frac{9}{2}, \infty}(\Rt)$ by H. Kozono \emph{et. al.} in \cite{KzonoAL}, but in this more general space  some \emph{supplementary hypothesis} were needed to obtain $\U=0$. Indeed, in Theorem $1.2$ of the article \cite{KzonoAL}  it is proven that if $\U \in L^{\frac{9}{2}, \infty}(\Rt)$ then we have the estimate $\int_{\Rt} \vert \vec{\nabla} \otimes \U (x)\vert^{2}dx \leq c\Vert \U \Vert^{3}_{L^{\frac{9}{2}, \infty}}$ and the desired identity $\U=0$ is then obtained under the  hypothesis 
\begin{equation}\label{cond-Kozono}
\Vert \U \Vert^{3}_{L^{\frac{9}{2}, \infty}} \leq \delta \int_{\Rt} \vert \vec{\nabla} \otimes \U (x)\vert^{2}dx,
\end{equation}  with $\delta >0$ small enough.  Although this supplementary hypothesis allow us to prove that $\U=0$ we may observe that it is a quite strong hypothesis and one of the aims of the article \cite{SerWang} by G. Seregin \& W. Wang is to relax the restriction imposed on the quantity $\Vert \U \Vert_{L^{\frac{9}{2}, \infty}}$. For this purpose, in Theorem $1.1$ of the article \cite{SerWang} the following result is proven:  if $\U$ is a smooth solution of equations (\ref{N-S-f-nulle-intro}) and if for a parameter $3<r<+\infty$ we have
\begin{equation}\label{Cond-Seregin-Wang}
M(r)= \sup_{R>1} R^{\frac{2}{3}-\frac{3}{r}} \Vert \U \Vert_{L^{r, \infty} (\mathcal{C}(R/2, R))} <+\infty,
\end{equation}   then we get the estimate $\int_{\Rt} \vert \vec{\nabla} \otimes \U (x) \vert^2 dx \leq c M^{3}(r)$, and moreover, if for $\delta>0$ small enough we have the supplementary \emph{a priori} control  $$M^{3}(r) \leq \delta \int_{\Rt} \vert \vec{\nabla} \otimes \U (x) \vert^2 dx,$$ then we get $\U=0$. Remark that for the value $r=\frac{9}{2}$ the condition $M^{3}(9/2) \leq \delta \int_{\Rt} \vert \vec{\nabla} \otimes \U (x) \vert^2 dx$ can be regarded as a relaxation of the condition (\ref{cond-Kozono}) given in \cite{KzonoAL}. \\

The first purpose of this article  is to review these results on the Liouville problem for stationary Navier-Stokes equations in the setting of Lorentz spaces. More precisely, we will prove that if we consider a slight smaller space than $L^{9/2, \infty}(\Rt)$: the space $L^{9/2, q}(\Rt)$ with $9/2 \leq q <+\infty$, then the information $\U \in L^{9/2, q}(\Rt)$ is sufficient to  derive the identity $\U=0$ and we do not need  any additional control on the quantity $\Vert \U \Vert_{L^{9/2, q}}$ contrary to the result given in \cite{KzonoAL}. Moreover, we will see that the space $L^{9/2, q}(\Rt)$ seems to be a \emph{critical} space to obtain the uniqueness of trivial solution  in the sens that if we have the information  $\U \in L^{r,q}(\Rt)$ for the values $\frac{9}{2}<r\leq q<+\infty$ then a faster decay condition on the solution $\U$ is required to obtain $\U=0$. \\

Our methods are based on a local estimate on the quantity $\int_{B_{R/2}} \vert \vec{\nabla} \otimes \U\vert^2 dx$ and this approach allows us to consider more general spaces than the Lorentz spaces. Thus, the second purpose of this article is to study the identity $\U=0$ in a the  framework  of the Morrey spaces $\dot{M}^{p,r}(\Rt)$ with $3\leq p<r<+\infty$, generalizing in this way some recent results. \\ 

This article is organized as follows: in Section \ref{Results} we state all the results obtained. In Section \ref{sec:estim-local} we prove a local estimate on the quantity above  from which we will able to study the Liouville problem in the setting of Lorentz space and this  will be done in Section \ref{sec:Lorentz}. Finally, in Section \ref{sec:Morrey} we extend our study to the setting of Morrey spaces. 

\section{Statement of the results}\label{Results} 
Recall first the definition of Lorentz spaces. Let  $f: \Rt \longrightarrow \mathbb{R}$ be a measurable function,  the  distribution function $d_{f}(\alpha)$ is defined as $$\ds{ d_{f}(\alpha) = dx \left( \{ x \in \Rt: \vert f(x) \vert >\alpha \} \right) },$$  where $dx$ denotes de Lebesgue measure. 
By definition, for $1\leq r < +\infty$ and $1\leq q \leq +\infty$  the Lorentz space $L^{r,q}(\Rt)$ is the space of measurable functions $f: \Rt \longrightarrow \mathbb{R}$ such that  $$\Vert f \Vert_{L^{r,q}}<+\infty,$$ where 
\begin{equation*}
 \Vert f \Vert_{L^{r,q}}= \left\{ \begin{array}{ll}  \vspace{3mm}
\ds{r^{1/q} \left( \int_{0}^{+\infty} \left(\alpha\, d^{1/r}_{f}(\alpha) \right)^{q} \frac{d \alpha}{\alpha}  \right)^{1/q}}, & \text{if} \,\, q<+\infty, \\
\ds{\sup_{\alpha>0} \left\{  \,  \alpha\,  d^{1/r}_{f}(\alpha)\right\}}, & \text{if} \,\, q=+\infty.  \end{array} \right. 
\end{equation*}
This space is a homogeneous space of degree $-\frac{3}{r}$ and we have the continuous embedding $L^{r} (\Rt) = L^{r,r}(\Rt) \subset L^{r,q}(\Rt)$ for $r < q \leq +\infty$. \\

In the framework of Lorentz spaces our first result is stated as follows:    
\begin{Theoreme}\label{Theo:Lorentz} Let $\U \in L^{2}_{loc}(\Rt)$ be a weak solution of the stationary Navier-Stokes  equations (\ref{N-S-f-nulle-intro}). 
\begin{enumerate}
\item[1)] If $\U \in L^{9/2, q}(\Rt)$, with $9/2 \leq q <+\infty$, then we have $\U=0$. 
\item[2)] If $\U \in L^{r,q}$, with $9/2<r\leq q <+\infty$, and moreover, if
	\begin{equation}\label{HIP-1}
	\sup_{R>1} R^{2-\frac{9}{r}} \Vert \U \Vert_{L^{r,q}(\mathcal{C}(R/4, 2R))} <+\infty,
	\end{equation} 
	 then we have $\U=0$. 
\end{enumerate}		
\end{Theoreme}

Several remarks follow from this result. First, as mentioned in the introduction, the result given in point $1)$  is of particular interest  since  this result can be regarded as a  improvement  of the results given in \cite{KzonoAL} and \cite{SerWang}. Moreover, due to the embedding $L^{9/2}(\Rt) \subset L^{9/2,q}(\Rt)$, Galdi's result \cite{Galdi} follows from this theorem. \\

Now, in point $2)$ we may observe that for the values $\frac{9}{2}<r<+\infty$ the information $\U \in L^{r,q}(\Rt)$ seems to be not enough to  prove that $\U=0$ and then  it is necessary a faster decay of the solution which is given in expression (\ref{HIP-1}). In this expression we may observe that as long as the parameter $r$ is larger than the critical value  $\frac{9}{2}$ the solution must have a faster decaying  at infinity. \\

As pointed out in the introduction, we also generalize our results to the framework of Morrey spaces and  we start by recalling their definition. For $1<p<r<+\infty$ the homogeneous Morrey space $\dot{M}^{p,r}(\Rt)$ is the set of functions $f \in L^{p}_{loc}(\Rt)$ such that 
\begin{equation}\label{Def-Morrey}
 \Vert f \Vert_{\dot{M}^{r,p}}= \sup_{R>0,\,\, x_0 \in \Rt}  R^{\frac{3}{r}- \frac{3}{p}} \left(  \int_{B(x_0, R)} \vert f (x) \vert^p dx \right)^{\frac{1}{p}}< +\infty,
\end{equation} where $B(x_0,R)$ denotes the ball centered at $x_0$ and with radio $R$.  This is a homogeneous space of degree $-\frac{3}{r}$ and moreover we have the following chain of continuous embeddings $L^{r}(\Rt)\subset L^{r, q}(\Rt) \subset L^{r, \infty}(\Rt) \subset \dot{M}^{p,r}(\Rt)$.  In the framework of Morrey spaces  our second result is the following:
   
\begin{Theoreme}\label{Theo:Morrey-sous-critique}  Let $\U \in L^{2}_{loc}(\Rt)$ be a weak solution of the stationary Navier-Stokes  equations (\ref{N-S-f-nulle-intro}). If $\U \in \dot{M}^{p,r}(\Rt)$ with $3\leq p < r < \frac{9}{2}$, then $\U=0$. 
\end{Theoreme}	

Observe that this result contains as particular case the uniqueness of the trivial solution of equations (\ref{N-S-f-nulle-intro}) in the setting of Lebesgue spaces $L^{r}(\Rt)$ and Lorentz spaces $L^{r,\infty}(\Rt)$ with the values $3<r<\frac{9}{2}$, and this fact extend to a more general framework some recent results obtained in the article \cite{ChJaLem}.\\

Now, It is natural to ask what happens for the values $\frac{9}{2} \leq r <+\infty$. For those values of parameter $r$,   following some ideas of the articles \cite{KzonoAL} and \cite{SerWang} exposed in the introduction, in our third result we prove some estimates of the quantity $\int_{\Rt} \vert \vec{\nabla} \otimes \U (x) \vert^2 dx$ by means of the quantity $\Vert \U \Vert_{\dot{M}^{p,r}}$ (where $3 \leq p < r$ and $\frac{9}{2} \leq r < +\infty$), and thus, the information $\U \in \dot{M}^{p,r}$ allow us to derive the fact that $\U \in \dot{H}^{1}(\Rt)$.  \\

\begin{Theoreme}\label{Theo:Morrey} Let $\U \in L^{2}_{loc}(\Rt)$ be a weak solution of the stationary Navier-Stokes  equations (\ref{N-S-f-nulle-intro}).  Suppose that $\U \in \dot{M}^{p,r}(\Rt)$ with $3 \leq p<r$ and $\frac{9}{2}\leq r <+\infty$.  
\begin{enumerate}
\item[1)] For  the limit value $r= \frac{9}{2}$ we have $\ds{\int_{\Rt} \vert \vec{\nabla} \otimes \U (x)\vert^2 dx \leq c \Vert \U \Vert^{3}_{\dot{M}^{p,\frac{9}{2}}}}$.  
\item[2)] For the values $\frac{9}{2}  < r <+\infty$, if  moreover  
\begin{equation}\label{HIP-2}
N(r)=\sup_{R>1} R^{2- \frac{9}{r}} \left( R^{\frac{3}{r} -\frac{3}{p}}  \left( \int_{\mathcal{C}(R/2,R)} \vert \U (x) \vert^{p} dx \right)^{\frac{1}{p}}\right)<+\infty,
\end{equation}  then we have $\ds{ \int_{\Rt} \vert \vec{\nabla} \otimes \U(x) \vert^2 dx \leq c \Vert \U \Vert^{2}_{\dot{M}^{p,r}} N(r)}$.  	 
\end{enumerate}	
\end{Theoreme}

Comparing this result with the results obtained in \cite{KzonoAL} and \cite{SerWang} (in the setting of Lorentz spaces) we may observe that point $1)$ below  generalizes  to Morrey spaces of the result given in    \cite{KzonoAL}, whereas if we compare  the expression (\ref{Cond-Seregin-Wang}) with the expression  (\ref{HIP-2}) below then we may see that  point $2)$ is in a certain sens a  generalization  to Morrey spaces of the result given in   \cite{SerWang}.\\

Now, in order to obtain the desired identity $\U=0$ in the framework of this result, and to the best of our knowledge, it is still necessary to make supplementary hypothesis on the solution $\U$. Following always the ideas of   \cite{KzonoAL} and \cite{SerWang} we could suppose an additional control on the quantities $\Vert \U \Vert_{\dot{M}^{p,\frac{9}{2}}}$ and $N(r)$ by means of $\int_{\Rt} \vert \vec{\nabla} \otimes \U(x)\vert^2 dx$, however  we will use here a different approach. 

\begin{Corollaire}\label{Corollary-Morrey}  Within the framework of Theorem \ref{Theo:Morrey}. If $\U \in \dot{B}^{-1}_{\infty,\infty}(\Rt)$ then we have $\U=0$.
\end{Corollaire}

Recall that the Besov space $\dot{B}^{-1}_{\infty,\infty}(\Rt)$, which is characterized as the set of distributions $f\in \mathcal{S}'(\Rt)$ such that $\Vert f \Vert_{\dot{B}^{-1}_{\infty, \infty}}= \ds{ \sup_{t>0} t^{\frac{1}{2}} \Vert h_t \ast f \Vert_{L^{\infty}}<+\infty}$ and where $h_t$ denotes the heat kernel, plays a very important role in the analysis on the Navier-Stokes equations (stationary and non stationary) since this is the largest space which is invariant under scaling properties of these equations (see the article \cite{JBourgainAL} and the books \cite{PGLR1} and \cite{PGLR2} for more references). Thus, in order obtain the identity $\U=0$, we have supposed  $\U \in \dot{B}^{-1}_{\infty, \infty}(\Rt)$ which is a condition on $\U$ less restrictive  compared to those made in \cite{KzonoAL} and \cite{SerWang}.  
\section{A local estimate}\label{sec:estim-local}
From now on $\U\in L^{2}_{loc}(\Rt)$ will be a weak solution of the stationary Navier-Stokes equations (\ref{N-S-f-nulle-intro}). Our results deeply relies on the following technical estimate (also known as a Caccioppoli type inequality):
\begin{Proposition}\label{Prop-Base} If the solution $\U$ verifies $\U \in L^{p}_{loc}(\Rt)$ and $\vec{\nabla} \otimes \U \in L^{\frac{p}{2}}_{loc}(\Rt)$   with  $3 \leq p< +\infty$, then for all $R>1$ we have 
	\begin{equation}\label{Ineq-base}
	\begin{split}
	\int_{B_{R/2}} \vert \vec{\nabla} \otimes \U \vert^2 dx  \leq  c \left(  \left( \int_{\mathcal{C}(R/2,R)} \vert \vec{\nabla}\otimes \U\vert^{\frac{p}{2}}dx \right)^{\frac{2}{p}} + \left( \int_{\mathcal{C}(R/2,R)} \vert \U \otimes \U\vert^{\frac{p}{2}}dx \right)^{\frac{2}{p}} \right)    \\ 
	  \times R^{2-\frac{9}{p}} \left(  \int_{\mathcal{C}(R/2,R)} \vert \U \vert^{p} dx \right)^{\frac{1}{p}}. 
	  \end{split}
	\end{equation}
\end{Proposition}	
\pv 
We start by introducing the test functions $\varphi_R$ and $\W_R$ as follows:  for a fixed $R>1$, we define first the function $\varphi_R\in \mathcal{C}^{\infty}_{0}(\Rt)$ by $0\leq \varphi_R\leq 1$ such that  for $\vert x \vert \leq  \frac{R}{2}$ we have $\varphi_R(x)=1$,  for $\vert x \vert \geq R$ we have $\varphi_R(x)=0$, and 
\begin{equation}\label{control-norme-test}
\Vert \vec{\nabla} \varphi_R \Vert_{L^{\infty}} \leq \frac{c}{R}. 
\end{equation}

 Next we define the function $\W_R$ as the solution of the problem
\begin{equation}\label{eq_W_R}
div(\W_R)=\vec{\nabla}\varphi_R\cdot \U, \quad \text{over}\,\, B_R, \quad \text{and}\quad \W_R=0 \,\, \text{over}\,\,  \partial B_R \cup \partial B_{\frac{R}{2}},
\end{equation} where  $\partial B_R= \{ x\in \Rt: \vert x \vert =R\}$.
Existence of such function $\W_R$ is assured by Lemma  $III. 3.1$ (page 162) of the book  \cite{Galdi} and where it is proven that for $1<q<+\infty$ we have $\W_R\in W^{1,q}(B_R)$ with $supp\,(\W_R)\subset \mathcal{C}(R/2, R)$ (the function $\W_R$ is extended by zero outside the set $\mathcal{C}(R/2, R)$) and 
\begin{equation}\label{contril-Lq-test}
\Vert \vec{\nabla}\otimes \W_R\Vert_{L^q(\mathcal{C}(R/2, R))}\leq c \Vert   \vec{\nabla}\varphi_R\cdot \U \Vert_{L^q(\mathcal{C}(R/2, R))}.
\end{equation} 

Once we have defined the functions $\varphi_R$ and $\W_R$ above, we consider now the function $\varphi_R \U-\W_R$ and we write  
\begin{equation}\label{eq01}
 \int_{B_R} \left( -\Delta \U +(\U \cdot \vec{\nabla})\U +\vec{\nabla}P\right)\cdot \left( \varphi_R \U-\W_R \right)dx=0.
\end{equation}

Remark that as $\U \in L^{p}_{loc}(\Rt)$ (with $3\leq p < +\infty$) then $\U \in L^{3}_{loc}(\Rt)$ and  by Theorem X.1.1 (page 658) of the book \cite{Galdi} we have $\U \in \mathcal{C}^{\infty}(\Rt)$ and $P\in \mathcal{C}^{\infty}(\Rt)$,  thus every term in the last identity above is well defined.\\

In identity (\ref{eq01}), we start by studying the third term in the left-hand side 
 and integrating  by parts we obtain
$$ \int_{B_R} \vec{\nabla}P \cdot \left( \varphi_R \U-\W_R \right)dx=-\int_{B_R}P\left( \vec{\nabla}\varphi_R\cdot \U+\varphi_R div(\U)-div(\W_R)\right)dx,$$
but since $\W_R$ is a solution of problem (\ref{eq_W_R}) and since $div(\U)=0$ we can write  $$\displaystyle{\int_{B_R}\vec{\nabla}P\cdot \left( \varphi_R \U-\W_R \right)dx=0},$$ and thus identity (\ref{eq01}) can be written as 
\begin{equation}\label{eq02}
\int_{B_R} -\Delta \U \cdot \left( \varphi_R \U-\W_R \right)dx + \int_{B_R}\left((\U \cdot \vec{\nabla})\U \right)\cdot \left( \varphi_R \U-\W_R \right)dx=0. 
\end{equation}

In this identity we study now the first term in the left-hand side  and always by integration by parts we have
\begin{eqnarray*}
& & \int_{B_R} -\Delta \U \cdot \left(\varphi_R \U- \W_R \right) dx  =  \sum_{i,j=1}^{3} \int_{B_R} (\partial_j U_i ) \partial_j (\varphi_R U_i-(W_R)_i)dx \\
&=& \sum_{i,j=1}^{3} \int_{B_R} (\partial_j U_i) (\partial_j \varphi_R) U_i dx +  \sum_{i,j=1}^{3}\int_{B_R}  \varphi_R (\partial_j U_i)^2 dx - \sum_{i,j=1}^{3} \int_{B_R} (\partial_j U_i ) \partial_j(W_R)_i dx \\
&=& \sum_{i,j=1}^{3}\int_{B_R} (\partial_j U_i) (\partial_j \varphi_R) U_i dx  + \int_{B_R} \varphi_R \vert \vec{\nabla} \otimes \U \vert^{2} dx - \sum_{i,j=1}^{3} \int_{B_R} (\partial_j U_i ) \partial_j(W_R)_i dx. 
\end{eqnarray*}

With this identity at hand, we get back to equation (\ref{eq02}) and we can write 
\begin{eqnarray*}
& & \sum_{i,j=1}^{3} \int_{B_R}(\partial_j U_i) (\partial_j \varphi_R) U_i dx + \int_{B_R} \varphi_R \vert \vec{\nabla} \otimes \U \vert^{2} dx - \sum_{i,j=1}^{3} \int_{B_R} (\partial_j U_i ) \partial_j(W_R)_i dx \\
&+ & \int_{B_R}\left((\U \cdot \vec{\nabla})\U \right)\cdot \left( \varphi_R \U-\W_R \right)dx=0,
\end{eqnarray*} hence we have 
\begin{eqnarray*}
\int_{B_R} \varphi_R \vert \vec{\nabla} \otimes \U \vert^{2} dx &=&  - \sum_{i,j=1}^{3}\int_{B_R} (\partial_j U_i) (\partial_j \varphi_R) U_i dx + \sum_{i,j=1}^{3} \int_{B_R} (\partial_j U_i ) \partial_j(W_R)_i dx \\
& & + \int_{B_R}\left((\U \cdot \vec{\nabla})\U \right)\cdot \left( \varphi_R \U-\W_R \right)dx. 
\end{eqnarray*}

But recall  the fact that test function $\varphi_R$ verifies $\varphi_R(x)=1$ if $\vert x \vert <\frac{R}{2}$, and  then we have 
$$ \int_{B_{R/2}} \vert \vec{\nabla} \otimes \U \vert^{2} dx \leq \int_{B_R} \varphi_R \vert \vec{\nabla} \otimes \U \vert^{2} dx.$$ Thus by this inequality and the identity above we can write the following estimate: 
\begin{eqnarray} \label{eq03} \nonumber
 \int_{B_{R/2}} \vert \vec{\nabla} \otimes \U \vert^{2} dx & \leq &  	- \sum_{i,j=1}^{3}\int_{B_R} (\partial_j U_i) (\partial_j \varphi_R) U_i dx + \sum_{i,j=1}^{3} \int_{B_R} (\partial_j U_i ) \partial_j(W_R)_i dx  \\ \nonumber
 & & 
 + \int_{B_R}\left((\U \cdot \vec{\nabla})\U \right)\cdot \left( \varphi_R \U-\W_R \right)dx \\
 & =& I_1+ I_2 +I_3. 
\end{eqnarray}

We study now these three terms above. In term $I_1$ remark that we have the function $\partial_i \varphi_R$, but  since the test function $\varphi_R$ verifies $\varphi_R(1)$ if $\vert x \vert <\frac{R}{2}$ and $\varphi_R(x)=0$  if $\vert x \vert > R$ then we have $supp\,(\vec{\nabla} \varphi_R) \subset \mathcal{C}(R/2,R)$, and thus we can write 
$$ I_1= -\sum_{i,j=1}^{3} \int_{\mathcal{C}(R/2,R)} \partial_j U_i (\partial_j \varphi_R) U_i dx.$$ 

Then, applying the Hold\"er inequalities with the relation $1=\frac{2}{p}+\frac{1}{q}$ we write 
\begin{eqnarray}\label{eq06} \nonumber 
 I_1 &\leq & \sum_{i,j=1}^{3} \left( \int_{\mathcal{C}(R/2,R)} \vert \partial_j U_i \vert^{\frac{p}{2}} dx \right)^{\frac{2}{p}} \left( \int_{\mathcal{C}(R/2,R)} \vert (\partial_j \varphi_R) U_i \vert^{q} dx \right)^{\frac{1}{q}}\\ \nonumber 
 &\leq  &c  \left( \int_{\mathcal{C}(R/2,R)}\vert \vec{\nabla} \otimes \U \vert^{\frac{p}{2}} dx \right)^{\frac{2}{p}} \Vert  \vec{\nabla} \varphi_R \Vert_{L^{\infty}} \left( \int_{\mathcal{C}(R/2,R)} \vert \U \vert^{q} dx\right)^{\frac{1}{q}} \\
 &\leq & c \left( \int_{\mathcal{C}(R/2,R)}\vert \vec{\nabla} \otimes \U \vert^{\frac{p}{2}} dx \right)^{\frac{2}{p}}  \underbrace{\frac{1}{R} \left( \int_{\mathcal{C}(R/2,R)} \vert \U \vert^{q} dx\right)^{\frac{1}{q}}}_{(a)} ,
\end{eqnarray} where the last estimate  is due to  (\ref{control-norme-test}).  We need to study now the term (a). Remark the fact that as $3 \leq p < +\infty$ and by the relation $1= \frac{2}{p}+ \frac{1}{q}$ then we have $q \leq 3 \leq p$, and thus we can write 
\begin{equation} \label{eq04} 
(a) \leq  \frac{R^{3(\frac{1}{q}-\frac{1}{p})}}{R}   \left( \int_{\mathcal{C}(R/2,R)} \vert \U \vert^{p} dx \right)^{\frac{1}{p}} \leq \frac{R^{3((
		1-\frac{2}{p}) - \frac{1}{p})}}{R}   \left( \int_{\mathcal{C}(R/2,R)} \vert \U \vert^{p} dx \right)^{\frac{1}{p}}  \leq R^{2-\frac{9}{p}} \left(  \int_{\mathcal{C}(R/2,R)} \vert \U \vert^{p} dx \right)^{\frac{1}{p}}. 
\end{equation}

With this estimate at hand  we write 
\begin{equation}\label{eq07}
I_1  \leq  c \left( \int_{\mathcal{C}(R/2,R)}\vert \vec{\nabla} \otimes \U \vert^{\frac{p}{2}} dx \right)^{\frac{2}{p}}  R^{2-\frac{9}{p}} \left(  \int_{\mathcal{C}(R/2,R)} \vert \U \vert^{p} dx \right)^{\frac{1}{p}}
\end{equation} 

In order to study  the term $I_2$ in (\ref{eq03}), recall that the have $supp\,(\W_R) \subset \mathcal{C}(R/2,2)$, hence we get $ supp \, (\vec{\nabla} \otimes \W_R) \subset  \mathcal{C}(R/2,R)$ and then we can write 
$$ I_2= \sum_{i,j=1}^{3} \int_{B_R} (\partial_j U_i ) \partial_j(W_R)_i dx =\sum_{i,j=1}^{3} \int_{\mathcal{C}(R/2,R)} (\partial_j U_i ) \partial_j(W_R)_i dx.$$

Now,  we apply the H\"older inequalities always with the relation $1= \frac{2}{p}+ \frac{1}{q}$ and we write
\begin{equation*}
I_2  \leq  c \left( \int_{\mathcal{C}(R/2,R)} \vert \vec{\nabla} \otimes \U \vert^{\frac{p}{2}} dx \right)^{\frac{2}{p}} \left( \int_{\mathcal{C}(R/2,R)} \vert \vec{\nabla} \otimes \W_R \vert^{q}dx \right)^{\frac{1}{q}},
\end{equation*} where it remains to study the second  term in the right-hand.  For this, applying first the estimate (\ref{contril-Lq-test}), then  applying the estimate (\ref{control-norme-test}) and finally by estimate (\ref{eq04}) we can write  
\begin{eqnarray*}
\left( \int_{\mathcal{C}(R/2,R)} \vert \vec{\nabla} \otimes \W_R \vert^{q}dx \right)^{\frac{1}{q}} & \leq & c\left( \int_{\mathcal{C}(R/2,R)} \vert \vec{\nabla} \varphi_R \cdot \U \vert^{q} dx \right)^{\frac{1}{q}} \leq c \Vert \vec{\nabla} \varphi_R \Vert_{L^{\infty}} \left( \int_{\mathcal{C}(R/2,R)} \vert  \U \vert^{q} dx \right)^{\frac{1}{q}} \\
&\leq & c\, R^{2-\frac{9}{p}} \left(  \int_{\mathcal{C}(R/2,R)} \vert \U \vert^{p} dx \right)^{\frac{1}{p}},
\end{eqnarray*}
and thus we have 
\begin{equation}\label{eq08}
I_2 \leq c \left( \int_{\mathcal{C}(R/2,R)} \vert \vec{\nabla} \otimes \U \vert^{\frac{p}{2}} dx \right)^{\frac{2}{p}}  R^{2-\frac{9}{p}} \left(  \int_{\mathcal{C}(R/2,R)} \vert \U \vert^{p} dx \right)^{\frac{1}{p}}. 
\end{equation} 

Finally we study the term $I_3$ in (\ref{eq03}). As we have $div(\U)=0$ then in this term we write $(\vec{\nabla} \cdot \U) \U= div(\U \otimes \U)$ and  we obtain
\begin{equation*}
I_3= \sum_{i,j=1}^{3} \int_{B_R} \partial_j (U_i U_j) (\varphi_R U_i -(W_R)_i) dx, 
\end{equation*} integrating by parts we write 
\begin{eqnarray}\label{eq09} \nonumber 
I_3& =& -\sum_{i,j=1}^{3} \int_{B_R} (U_i U_j) ((\partial_j \varphi_R) U_i + \varphi_R (\partial_j U_i) -\partial_j(W_R)_i) dx \\ \nonumber
&=& - \sum_{i,j=1}^{3}  \int_{B_R} (U_i U_j) ((\partial_j \varphi_R) U_i dx -\sum_{i,j=1}^{3} \int_{B_R}  (U_i U_j) \varphi_R (\partial_j U_i) dx + \sum_{i,j=1}^{3} \int_{B_R} (U_i U_j) \partial_j(W_R)_i) dx\\
&=& I_{3,a}+ I_{3,b}+I_{3,c},
\end{eqnarray}  where  we will study these three terms separately. In term $I_{3,a}$, as we have $supp\,(\vec{\nabla} \varphi_R) \subset  \mathcal{C}(R/2,R)$ then we write 
$$ I_{3,a}=-  \sum_{i,j=1}^{3}  \int_{\mathcal{C}(R/2,R)} (U_i U_j) ((\partial_j \varphi_R) U_i dx,$$ then, applying first  the H\"older inequalities (with the same relation $1= \frac{2}{p}+\frac{1}{q}$) and thereafter, applying first estimate (\ref{eq06}) and then estimate (\ref{eq04})   we have 
\begin{eqnarray}\label{eq10}\nonumber
 I_{3,a} &\leq & c \left( \int_{\mathcal{C}(R/2,R)} \vert \U\otimes \U\vert^{\frac{p}{2}}dx \right)^{\frac{2}{p}} \left( \int_{\mathcal{C}(R/2,R)} \vert \vec{\nabla} \varphi_R \cdot \U \vert^{q} dx\right)^{\frac{1}{q}}  \\
 &\leq &c \left( \int_{\mathcal{C}(R/2,R)} \vert \U\otimes \U\vert^{\frac{p}{2}}dx \right)^{\frac{2}{p}}  R^{2-\frac{9}{p}} \left(  \int_{\mathcal{C}(R/2,R)} \vert \U \vert^{p} dx \right)^{\frac{1}{p}}.
 \end{eqnarray}
 
In order to estimate the term $I_{3,b}$ we write 
\begin{equation*}
I_{3,b} = -\sum_{i,j=1}^{3} \int_{B_R}  (U_i U_j) \varphi_R (\partial_j U_i) dx= -\sum_{i,j=1}^{3} \int_{B_R}  U_j \varphi_R ( (\partial_j U_i) U_i)  dx = - \frac{1}{2} \sum_{i,j=1}^{3} \int_{B_R}  U_j \varphi_R \partial_j (U^{2}_{i}) dx, 
\end{equation*} then, by integration by parts, and moreover, using the fact that $div(\U)=0$ and since the function $\vec{\nabla} \varphi_R$ is localized at the set $\mathcal{C}(R/2, R)$, then we get: 
\begin{eqnarray*}
& & - \frac{1}{2} \sum_{i,j=1}^{3} \int_{B_R}  U_j \varphi_R \partial_j (U^{2}_{i}) dx= \frac{1}{2}	 \sum_{i,j=1}^{3} \int_{B_R} \partial_j(U_j \varphi_R) U^{2}_{i} dx = \frac{1}{2}   \sum_{i,j=1}^{3} \int_{B_R} (\partial_j U_j) \varphi_R U^{2}_{i} dx \\
& & + \frac{1}{2}   \sum_{i,j=1}^{3} \int_{B_R} U_j (\partial_j \varphi_R ) U^{2}_{i} dx = \frac{1}{2} \int_{B_R} div(\U) \varphi_R \vert \U \vert^2 dx + \frac{1}{2} \sum_{i,j=1}^{3} \int_{\mathcal{C}(R/2,R)} U_j (\partial_j \varphi_R) U^{2}_{i} dx\\
&=& \frac{1}{2} \sum_{i,j=1}^{3} \int_{\mathcal{C}(R/2,R)} U_j (\partial_j \varphi_R) U^{2}_{i} dx= \frac{1}{2} \sum_{i,j=1}^{3} \int_{\mathcal{C}(R/2,R)} U^{2}_{i} (\partial_j \varphi_R) U_j dx. 
\end{eqnarray*}

With this identity at hand and following the same estimates done for the term $I_{3,a}$ in (\ref{eq10}) we have   
\begin{equation}\label{eq11}
I_{3,b}  \leq  c \left( \int_{\mathcal{C}(R/2,R)} \vert \U\otimes \U\vert^{\frac{p}{2}}dx \right)^{\frac{2}{p}}  R^{2-\frac{9}{p}} \left(  \int_{\mathcal{C}(R/2,R)} \vert \U \vert^{p} dx \right)^{\frac{1}{p}}. 
\end{equation}

Now, in order to study term $I_{3,c}$ remark that  using the inequality (\ref{contril-Lq-test}) and following always the estimates done for term $I_{3,a}$ (see (\ref{eq10})) we have 
\begin{equation}\label{eq20}
I_{3,c} \leq c  \left( \int_{\mathcal{C}(R/2,R)} \vert \U\otimes \U\vert^{\frac{p}{2}}dx \right)^{\frac{2}{p}}  R^{2-\frac{9}{p}} \left(  \int_{\mathcal{C}(R/2,R)} \vert \U \vert^{p} dx \right)^{\frac{1}{p}}.
\end{equation}

With estimates (\ref{eq10}),  (\ref{eq11}) and (\ref{eq20}) we get back to identity (\ref{eq09}) hence  we have 
\begin{equation}\label{eq12}
I_3 \leq c \left( \int_{\mathcal{C}(R/2,R)} \vert \U\otimes \U\vert^{\frac{p}{2}}dx \right)^{\frac{2}{p}}  R^{2-\frac{9}{p}} \left(  \int_{\mathcal{C}(R/2,R)} \vert \U \vert^{p} dx \right)^{\frac{1}{p}}.
\end{equation} 

Finally, once we dispose of estimates (\ref{eq07}), (\ref{eq08}) and (\ref{eq12}), applying these estimates in each term in the right-hand side of  (\ref{eq03})  we obtain the desired estimate (\ref{Ineq-base}).  \finpv
\section{The Lorentz spaces: proof of Theorem \ref{Theo:Lorentz}}\label{sec:Lorentz}
Suppose that the solution $\U \in L^{2}_{loc}(\Rt)$ of equations (\ref{N-S-f-nulle-intro}) verifies $\U \in L^{r,q}(\Rt)$ with $\frac{9}{2}\leq   r \leq q <+\infty$. The first think to do is to prove that  $\U$ verifies the hypothesis of Proposition \ref{Prop-Base}, and for this recall the following estimate: for $1<p<r \leq q < +\infty$ and for $R>1$ we have 
\begin{equation}\label{estim-Lorentz}
\int_{B_R} \vert \U \vert^{p} dx \leq c \, \, R^{3(1-\frac{p}{r})} \Vert \U \Vert^{p}_{L^{r,\infty}} \leq c   \, R^{3(1-\frac{p}{r})} \Vert \U \Vert^{p}_{L^{r,q}},
\end{equation}  see Proposition $1.1.10$, page 22 of the  book \cite{DCh} for a proof of this fact. From this estimate we have $\U \in L^{p}_{loc}(\Rt)$ and then it remains to prove that $\vec{\nabla} \otimes \U \in L^{\frac{p}{2}}_{loc}(\Rt)$ for $3\leq p <+\infty$.  Indeed, since $\U$ verifies the equations (\ref{N-S-f-nulle-intro}) and since $div(\U)=0$ then this solution  can be written as follows 
$$ \U = - \frac{1}{\Delta} \left( \P \left(( \U \cdot \vec{\nabla}) \U \right) \right)=\sum_{j=1}^{3} - \frac{1}{\Delta} \left( \P \left( \partial_j (U_j \U) \right) \right) ,$$ where $\P$ is the Leray's projector. Then, for $i=1,2,3$ we have 
\begin{equation}\label{eq16}
 \partial_i \U = -  \sum_{j=1}^{3}  \frac{1}{\Delta} \left( \P \left( \partial_i  \partial_j (U_j \U) \right) \right)= \sum_{j=1}^{3} \P \left( \mathcal{R}_{i} \mathcal{R}_{j} (U_j \U)  \right),
\end{equation}  where recall that $\mathcal{R}_{i}= \frac{\partial_i}{\sqrt{-\Delta}}$ denotes the i-th Riesz transform. Thus, by continuity of the operator $\P(\mathcal{R}_{i} \mathcal{R}_{j})$ on Lorentz spaces $L^{r,q}(\Rt)$ for the values $1< r \leq q$ (see the article \cite{Brandolese}) and applying the H\"older inequalities we obtain the following estimate: 
$$ \Vert \vec{\nabla} \otimes \U \Vert_{L^{\frac{r}{2}, \frac{q}{2}}} \leq c \sum_{i,j=1}^{3} \Vert \P(\mathcal{R}_i \mathcal{R}_j (U_j \U)) \Vert_{L^{\frac{r}{2}, \frac{q}{2}}} \leq c \Vert \U \otimes \U \Vert_{L^{\frac{r}{2},\frac{q}{2}}} \leq c \Vert \U \Vert^{2}_{L^{r,q}} .$$  With this estimate at hand  we can use now the last estimate in (\ref{estim-Lorentz}) (with $1<\frac{p}{2}<\frac{r}{2} <+\infty$) to write 
\begin{equation}\label{estim-Lorentz-2}
\int_{B_R} \vert \vec{\nabla} \otimes U \vert^{\frac{p}{2}} dx \leq c\, R^{3(1-\frac{p/2}{r/2})} \Vert \vec{\nabla} \otimes  \U \Vert^{\frac{p}{2}}_{L^{\frac{r}{2},\frac{q}{2}}}, 
\end{equation} hence we obtain $\vec{\nabla} \otimes \U \in L^{\frac{p}{2}}_{loc}(\Rt)$. \\

Thus,  by  Proposition \ref{Prop-Base} the solution $\U$ verifies (\ref{Ineq-base}) and by this estimate  we can write for all $R>1$ 
\begin{eqnarray} \nonumber
& & \int_{B_{\frac{R}{2}}} \vert \vec{\nabla} \otimes \U \vert^2 dx  \leq    \left(  \left( \int_{\mathcal{C}(R/2,R)} \vert \vec{\nabla}\otimes \U\vert^{\frac{p}{2}}dx \right)^{\frac{2}{p}} + \left( \int_{\mathcal{C}(R/2,R)} \vert \U \otimes \U\vert^{\frac{p}{2}}dx \right)^{\frac{2}{p}} \right)  R^{2-\frac{9}{p}} \left(  \int_{\mathcal{C}(R/2,R)} \vert \U \vert^{p} dx \right)^{\frac{1}{p}} \\  \label{eq17} 
&\leq & c \left(  \left( \frac{1}{R^3} \int_{\mathcal{C}(R/2,R)} \vert \vec{\nabla}\otimes \U\vert^{\frac{p}{2}}dx \right)^{\frac{2}{p}} + \left(\frac{1}{R^3} \int_{\mathcal{C}(R/2,R)} \vert \U \otimes \U\vert^{\frac{p}{2}}dx \right)^{\frac{2}{p}} \right)  R^{2} \left( \frac{1}{R^3}  \int_{\mathcal{C}(R/2,R)} \vert \U \vert^{p} dx \right)^{\frac{1}{p}},
\end{eqnarray} hence we have 
\begin{equation} \label{eq13}
\begin{split}
\int_{B_{\frac{R}{2}}} \vert \vec{\nabla} \otimes \U \vert^2 dx  \leq  c \underbrace{\left( R^{\frac{6}{r}} \left( \frac{1}{R^3} \int_{\mathcal{C}(R/2,R)} \vert \vec{\nabla}\otimes \U\vert^{\frac{p}{2}}dx \right)^{\frac{2}{p}} + R^{\frac{6}{r}} \left(\frac{1}{R^3} \int_{\mathcal{C}(R/2,R)} \vert \U \otimes \U\vert^{\frac{p}{2}}dx \right)^{\frac{2}{p}} \right)}_{(a)}  \\
 \times  R^{2-\frac{9}{r}}  \underbrace{\left( R^{\frac{3}{r}} \left( \frac{1}{R^3}  \int_{\mathcal{C}(R/2,R)} \vert \U \vert^{p} dx \right)^{\frac{1}{p}}\right)}_{(b)}, 
 \end{split}
\end{equation} where we will estimate the terms $(a)$ and $(b)$. For  this    we introduce  the cut-off function  $\theta_R \in \mathcal{C}^{\infty}_{0}(\Rt)$ such that  $\theta_R =1$ on  $ \mathcal{C}(R/2,R)$, $supp\, (\theta_R) \subset \mathcal{C}(R/4, 2R)$ and $\Vert \vec{\nabla} \theta_R \Vert_{L^{\infty}} \leq \frac{c}{R}$; and we consider the localized functions $\theta_R \U$ and $\theta_R (\vec{\nabla} \otimes \U)$. \\ 

Now, as we have $\theta_R=1$  on the set $ \mathcal{C}(R/2,R)$ then for the first term in $(a)$ we  can write 

$$\int_{\mathcal{C}(R/2, R)} \vert \vec{\nabla} \otimes \U \vert^{\frac{p}{2}}dx  = \int_{\mathcal{C}(R/2, R)}  \vert \theta_R( \vec{\nabla} \otimes \U) \vert^{\frac{p}{2}}dx  \leq \int_{B_{2R}}   \vert  \theta_R (\vec{\nabla} \otimes \U) \vert^{\frac{p}{2}}dx, $$

and applying the estimate (\ref{estim-Lorentz-2}) with the function $\theta_R (\vec{\nabla} \otimes \U)$ (and with $1<\frac{p}{2}<\frac{r}{2}<+\infty$) we have
$$   \int_{B_{2R}} \vert \theta_R (\vec{\nabla} \otimes \U)\vert^{\frac{p}{2}}dx \leq  c \, R^{3(1-\frac{p/2}{r/2})} \Vert \theta_R (\vec{\nabla} \otimes \U) \Vert^{\frac{p}{2}}_{L^{\frac{r}{2},\frac{q}{2}}} ,$$ hence the first term in expression (a) is estimated as   
$$ R^{\frac{6}{r}} \left( \frac{1}{R^3}\int_{\mathcal{C}(R/2,R)} \vert \vec{\nabla} \otimes \U \vert^{\frac{p}{2}}  dx \right)^{\frac{2}{p}} \leq  c \Vert \theta_R (\vec{\nabla} \otimes \U)\Vert_{L^{\frac{r}{2},\frac{q}{2}}}.$$ 

The second  term in $(a)$ treated in a similar way: first we write 
$$ \int_{\mathcal{C}(R/2,R)} \vert \U \otimes \U \vert^{\frac{p}{2}}dx = \int_{\mathcal{C}(R/2,R)} \vert (\theta_R \U)  \otimes (\theta_R \U) \vert^{\frac{p}{2}}dx  \leq \int_{B_{2R}} \vert (\theta_R \U) \otimes (\theta_R \U) \vert^{\frac{p}{2}} dx,$$ then we apply estimate (\ref{estim-Lorentz-2}) with the function $(\theta_R \U) \otimes (\theta_R \U) $ (always with $1<\frac{p}{2}<\frac{r}{2}<+\infty$) and by  the H\"older inequalities we have 
$$ \int_{B_{2R}} \vert (\theta_R \U) \otimes (\theta_R \U) \vert^{\frac{p}{2}} dx  \leq c  R^{3(1- \frac{p/2}{r/2})} \Vert (\theta_R \U) \otimes (\theta_R \U) \Vert^{\frac{p}{2}}_{L^{\frac{r}{2},\frac{q}{2}}} \leq c R^{3(1- \frac{p/2}{r/2})} \Vert \theta_R \U \Vert^{p}_{L^{r,q}},$$ hence we can write 
$$ R^{\frac{6}{r}} \left( \frac{1}{R^3}\int_{\mathcal{C}(R/2,R)} \vert \U \otimes \U \vert^{\frac{p}{2}}  dx \right)^{\frac{2}{p}} \leq  c \Vert \theta_R \U \Vert^{2}_{L^{r,q}}. $$

With these inequalities, the term $(a)$ above is estimated as follows:
\begin{equation}\label{eq14}
(a)\leq c  \left(\Vert \theta_R(\vec{\nabla} \otimes \U) \Vert_{L^{\frac{r}{2},\frac{q}{2}}} + \Vert \theta_R \U \Vert^{2}_{L^{r,q}}\right).
\end{equation}

We study now the term $(b)$. Following similar estimates done for  term $(a)$: applying always estimate (\ref{estim-Lorentz})  and as  $supp\, (\theta_R) \subset \mathcal{C}(R/4, 2R)$,  we can write 
$$ \int_{\mathcal{C}(R/2,R)} \vert \U \vert^p dx \leq c \int_{B_{2R}} \vert \theta_R \U \vert^p dx \leq c R^{3(1-\frac{p}{r})}\Vert \theta_R \U \Vert^{p}_{L^{r, q}} \leq c R^{3(1-\frac{p}{r})} \Vert \U \Vert^{p}_{L^{r,q}(\mathcal{C}(R/4,R))},$$ hence we obtain 
\begin{equation}\label{eq15}
(b) \leq   R^{\frac{3}{r}} \left(\frac{1}{R^3}  \int_{\mathcal{C}(R/2,R)} \vert \U \vert^{p} dx \right)^{\frac{1}{p}}  \leq c  \Vert \U \Vert_{L^{r,q}(\mathcal{C}(R/4,R))}. 
\end{equation}

Once we dispose of estimates (\ref{eq14}) and (\ref{eq15}) we get back to (\ref{eq13}) and we write 
\begin{equation}\label{eq18}
\int_{B_{\frac{B}{2}}} \vert \vec{\nabla} \otimes \U \vert^2 dx \leq c \left( \Vert \theta_R (\vec{\nabla} \otimes \U)\Vert_{L^{\frac{r}{2},\frac{q}{2}}} +  \Vert \theta_R \U \Vert^{2}_{L^{r,q}}\right) \, R^{2-\frac{9}{r}}    \Vert \U \Vert_{L^{r,q}(\mathcal{C}(R/4,R))},
\end{equation} and at this point we will consider two cases for the value of the parameters $r$ and $q$:
\begin{enumerate}
	\item[1)] For $ r=  \frac{9}{2}$ and $9/2 \leq q <+\infty$.  By   (\ref{eq18}) we can write 
	$$ \int_{B_{\frac{B}{2}}} \vert \vec{\nabla} \otimes \U \vert^2 dx \leq c \left( \Vert \theta_R (\vec{\nabla} \otimes \U)\Vert_{L^{\frac{9}{4},\frac{q}{2}}} +  \Vert \theta_R \U \Vert^{2}_{L^{9/2,q}}\right)    \Vert \U \Vert_{L^{\frac{9}{2},q}}.$$ 
	
	Now, recall that we have $supp\, (\theta_R) \subset \mathcal{C}(R/4,2R)$ and then we obtain $\ds{\lim_{R\longrightarrow +\infty} \theta_R \U =0}$ \emph{a.e.} in $\Rt$ and since we have  $q<+\infty$ we can apply   the dominated convergence theorem in Lorentz spaces (see Theorem $1.2.8$, page 74 of the book \cite{DCh}) to obtain $\ds{\lim_{R\longrightarrow +\infty} \Vert \theta_R (\vec{\nabla} \otimes \U) \Vert_{L^{9/4,q/2}}=0}$ and $\ds{\lim_{R\longrightarrow +\infty} \Vert \theta_R \U \Vert_{L^{9/2, q}}=0}$. Thus,  taking the limit $\ds{\lim_{R\longrightarrow +\infty}}$ in the estimate above we obtain $\Vert \vec{\nabla} \otimes \U \Vert^{2}_{L^2}=0$. Moreover, by the Hardy-Littlewood-Sobolev we also have $\Vert \U \Vert_{L^6} \leq c \Vert \vec{\nabla} \otimes \U \Vert_{L^2}$, hence we get the desired identity $\U=0$.
	\item[2)] For $\frac{9}{2}<r \leq q <+\infty$. In this case by estimate (\ref{eq18}) we have
	$$  \int_{B_{\frac{B}{2}}} \vert \vec{\nabla} \otimes \U \vert^2 dx \leq c \left( \Vert \theta_R (\vec{\nabla} \otimes \U)\Vert_{L^{\frac{r}{2},\frac{q}{2}}} +  \Vert \theta_R \U \Vert^{2}_{L^{r,q}}\right) \left( \sup_{R>1}   R^{2-\frac{9}{r}}    \Vert \U \Vert_{L^{r,qs}(\mathcal{C}(R/4,2R))}\right),$$ where by formula (\ref{HIP-1}) we know that the last  term in the right-hand side is bounded. Thus, always by the fact that $\ds{\lim_{R\longrightarrow +\infty} \Vert \theta_R (\vec{\nabla} \otimes \U) \Vert_{L^{\frac{r}{2},\frac{q}{2}}}=0}$ and $\ds{\lim_{R\longrightarrow +\infty} \Vert \theta_R \U \Vert_{L^{r,q}}=0}$ and taking the limit  $\ds{\lim_{R\longrightarrow +\infty}}$ in this estimate we obtain the identity $\U=0$.  Theorem \ref{Theo:Lorentz} is proven. \finpv  
\end{enumerate}	
\section{The Morrey spaces}\label{sec:Morrey}
Suppose now  the solution $\U \in L^{2}_{loc}(\Rt)$ of equations (\ref{N-S-f-nulle-intro}) verifies $\U \in \dot{M}^{p,r}(\Rt)$ with $3\leq p<r<+\infty$. Before to prove our results  we need to verify that the solution $\U$ satisfies the hypothesis of Proposition \ref{Prop-Base}:  $\U \in L^{p}_{loc}(\Rt)$  and $\vec{\nabla} \otimes \U \in L^{\frac{p}{2}}_{loc}(\Rt)$ with $3\leq p <+\infty$, but, as $\U \in \dot{M}^{p,r}(\Rt)$ then we have  
$\U \in L^{p}_{loc}(\Rt)$ (see Definition (\ref{Def-Morrey}) of Morrey spaces)  so it remains to verify that $\vec{\nabla} \otimes \U \in L^{\frac{p}{2}}_{loc}(\Rt)$ and  for this we will prove that $\vec{\nabla} \otimes \U \in \dot{M}^{\frac{p}{2}, \frac{r}{2}}(\Rt)$. Indeed, by identity (\ref{eq16}), the continuity of the operator $\P(\mathcal{R}_{i} \mathcal{R}_{j})$ on the Morrey spaces $\dot{M}^{r,p}(\Rt)$ with the values $1 <p<r<+\infty$  (see Lemma $4.2$ of the article \cite{Kato})  and applying the H\"older inequalities we can write 
\begin{equation}\label{eq19}
 \Vert \vec{\nabla} \otimes \U \Vert_{\dot{M}^{\frac{p}{2}, \frac{r}{2}}} \leq c \sum_{i,j=1}^{3} \Vert \P(\mathcal{R}_i \mathcal{R}_j (U_j \U)) \Vert_{\dot{M}^{\frac{p}{2}, \frac{r}{2}}} \leq c \Vert \U \otimes \U \Vert_{\dot{M}^{\frac{p}{2},\frac{r}{2}}} \leq c \Vert \U \Vert^{2}_{\dot{M}^{p,r}}.
 \end{equation}
 
Once we have the information $\U \in L^{p}_{loc}(\Rt)$  and $\vec{\nabla} \otimes \U \in L^{\frac{p}{2}}_{loc}(\Rt)$, by Proposition \ref{Prop-Base} we dispose of the inequality (\ref{Ineq-base}) and with this estimate at hand we will consider the following  cases of the values of parameters  $p$ and $r$.
\subsection{Proof of Theorem \ref{Theo:Morrey-sous-critique}} 
We consider here the values $3 \leq p<r<\frac{9}{2}$. By estimate (\ref{Ineq-base}) and following the same computations done in estimate (\ref{eq17}) we can write 
\begin{equation*}
\begin{split}
\int_{B_{\frac{R}{2}}} \vert \vec{\nabla} \otimes \U \vert^2 dx  \leq c \underbrace{\left(  \left( \frac{1}{R^3} \int_{\mathcal{C}(R/2,R)} \vert \vec{\nabla}\otimes \U\vert^{\frac{p}{2}}dx \right)^{\frac{2}{p}} + \left(\frac{1}{R^3} \int_{\mathcal{C}(R/2,R)} \vert \U \otimes \U\vert^{\frac{p}{2}}dx \right)^{\frac{2}{p}} \right)}_{(a)} \\
\times  \underbrace{R^{2} \left( \frac{1}{R^3}  \int_{\mathcal{C}(R/2,R)} \vert \U \vert^{p} dx \right)^{\frac{1}{p}}}_{(b)},
\end{split} 
\end{equation*}  where will estimate the terms $(a)$ and $(b)$. In term $(a)$ remark that as $\U \in \dot{M}^{p,r}(\Rt)$ then by (\ref{eq19}) we have $\vec{\nabla} \otimes \U \in \dot{M}^{\frac{p}{2}, \frac{r}{2}}(\Rt)$ and moreover by the H\"older inequalities we have $\U \otimes \U \in \dot{M}^{\frac{p}{2}, \frac{r}{2}}(\Rt)$. Thus, by definition of Morrey spaces  (see (\ref{Def-Morrey})) we can write 
\begin{equation*}
(a)\leq c \left(   \Vert \vec{\nabla} \otimes \U \Vert_{\dot{M}^{\frac{p}{2}, \frac{r}{2}}} +  \Vert \U \otimes \U \Vert_{\dot{M}^{\frac{p}{2}, \frac{r}{2}}} \right) R^{-\frac{6}{r}}. 
\end{equation*}
moreover, always by the fact that  $\U \in \dot{M}^{p,r}(\Rt)$  for term $(b)$ we write 
\begin{equation*}
(b)\leq R^{2}  \left( \Vert \U \Vert_{\dot{M}^{p,r}} R^{-\frac{3}{r}}\right) \leq  \Vert \U \Vert_{\dot{M}^{p,r}} R^{2-\frac{3}{r}}, 
\end{equation*} 
and with this estimates on terms $(a)$ and $(b)$ we obtain
\begin{equation*}
\int_{B_{\frac{R}{2}}} \vert \vec{\nabla} \otimes \U \vert^2 dx  \leq c \left(   \Vert \vec{\nabla} \otimes \U \Vert_{\dot{M}^{\frac{p}{2}, \frac{r}{2}}} +  \Vert \U \otimes \U \Vert_{\dot{M}^{\frac{p}{2}, \frac{r}{2}}} \right) \Vert  \U \Vert_{\dot{M}^{p,r}} R^{2-\frac{9}{r}}.
\end{equation*} 

But  recall that we have $3<r<\frac{9}{2}$ hence we get $-1<2-\frac{9}{r}<0$ and then, taking the limit $\ds{\lim_{R\longrightarrow +\infty}}$
we have $\Vert \vec{\nabla} \otimes \U \Vert^{2}_{L^2}=0$ hence we obtain the identity $\U=0$. Theorem \ref{Theo:Morrey-sous-critique} is now proven. \finpv

\subsection{Proof of Theorem \ref{Theo:Morrey}} 
\begin{enumerate}
	\item[1)] For the values $3\leq p<\frac{9}{2}$ and $r=\frac{9}{2}$. In this case we have $\U \in \dot{M}^{p,\frac{9}{2}}$. Following the same computations done in estimate  (\ref{eq13}) and moreover, always  by definition of the Morrey spaces   given in (\ref{Def-Morrey})  we get  the uniform bound 
	\begin{eqnarray*}
		\int_{B_{\frac{R}{2}}} \vert \vec{\nabla} \otimes \U \vert^2 dx & \leq & c \left( R^{\frac{4}{3}} \left( \frac{1}{R^3} \int_{\mathcal{C}(R/2,R)} \vert \vec{\nabla}\otimes \U\vert^{\frac{p}{2}}dx \right)^{\frac{2}{p}} + R^{\frac{4}{3}} \left(\frac{1}{R^3} \int_{\mathcal{C}(R/2,R)} \vert \U \otimes \U\vert^{\frac{p}{2}}dx \right)^{\frac{2}{p}} \right) \\
		& & \times   \left( R^{\frac{2}{3}} \left( \frac{1}{R^3}  \int_{\mathcal{C}(R/2,R)} \vert \U \vert^{p} dx \right)^{\frac{1}{p}}\right) \\
		&\leq & c \left( \Vert \vec{\nabla} \otimes \U \Vert_{\dot{M}^{\frac{p}{2}, \frac{9}{4}}} + \Vert \U \otimes \U \Vert_{\dot{M}^{\frac{p}{2}, \frac{9}{4}}}\right) \Vert  \U \Vert_{\dot{M}^{p, \frac{9}{2}}} \leq c \Vert \U \Vert^{3}_{\dot{M}^{p,\frac{9}{2}}},	
	\end{eqnarray*}	and  taking the limit $\ds{\lim_{R\longrightarrow +\infty}}$ we obtain $\ds{\int_{\Rt} \vert \vec{\nabla} \otimes \U \vert^2 dx \leq c \Vert \U \Vert^{3}_{\dot{M}^{p,\frac{9}{2}}}}$.   
\item[2)] For the values $3\leq p \leq \frac{9}{2}$ and  $\frac{9}{2} < r<+\infty$. Always by estimate (\ref{eq13})  for  all $R>1$ we write 
\begin{equation*}
\begin{split}
\int_{B_{\frac{R}{2}}} \vert \vec{\nabla} \otimes \U \vert^2 dx  \leq  c \underbrace{\left( R^{\frac{6}{r}} \left( \frac{1}{R^3} \int_{\mathcal{C}(R/2,R)} \vert \vec{\nabla}\otimes \U\vert^{\frac{p}{2}}dx \right)^{\frac{2}{p}} + R^{\frac{6}{r}} \left(\frac{1}{R^3} \int_{\mathcal{C}(R/2,R)} \vert \U \otimes \U\vert^{\frac{p}{2}}dx \right)^{\frac{2}{p}} \right)}_{(a)} \\
 \times  \underbrace{R^{2-\frac{9}{r}}  \left( R^{\frac{3}{r}} \left( \frac{1}{R^3}  \int_{\mathcal{C}(R/2,R)} \vert \U \vert^{p} dx \right)^{\frac{1}{p}}\right)}_{(b)},
\end{split} 
\end{equation*}
where,  as we have  $\U \in \dot{M}^{p,r}(\Rt)$ then the term $(a)$ is uniformly bounded as follows:
$$ (a)\leq c ( \Vert \vec{\nabla} \otimes \U \Vert_{\dot{M}^{\frac{p}{2}, \frac{r}{2}}} + \Vert \U \otimes \U \Vert_{\dot{M}^{\frac{p}{2}, \frac{r}{2}}})  \leq c \Vert \U \Vert^{2}_{\dot{M}^{p,r}},$$ moreover,  the term $(b)$ is uniformly bounded as 
$$ (b) \leq R^{2-\frac{9}{r}}  \left( R^{\frac{3}{r}} \left( \frac{1}{R^3}  \int_{\mathcal{C}(R/2,R)} \vert \U \vert^{p} dx \right)^{\frac{1}{p}}\right) \leq N(r),$$ where the quantity $N(r)<+\infty$ is defined in formula (\ref{HIP-2}).\\

With these estimates we can write $\ds{ \int_{B_{\frac{R}{2}}} \vert \vec{\nabla} \otimes \U \vert^2 dx  \leq  c \Vert \U \Vert^{2}_{\dot{M}^{p,r}} N(r)}$, and taking the limit $\ds{\lim_{R\longrightarrow +\infty}}$ we obtain $\ds{\int_{\Rt} \vert \vec{\nabla} \otimes \U \vert^2 dx \leq c \Vert \U \Vert^{2}_{\dot{M}^{p,r}} N(r)}$. Theorem \ref{Theo:Morrey} is  proven. \finpv
\end{enumerate}
\subsection{Proof of Corollary \ref{Corollary-Morrey}} 

 As $\int_{\Rt} \vert \vec{\nabla} \otimes \U \vert^2 dx <+\infty$ we get  $\U \in \dot{H}^{1}(\Rt)$, and with the information  $\U \in \dot{B}^{-1}_{\infty, \infty}(\Rt)$ we can apply the improved Sobolev inequalities (see the article \cite{GerardMeyerOru} for a proof of these inequalities) and we write $\ds{\Vert \U \Vert_{L^4} \leq c \Vert \U \Vert^{\frac{1}{2}}_{\dot{H}^{1}} \Vert \U \Vert^{\frac{1}{2} }_{\dot{B}^{-1}_{\infty,\infty}}}$. Once we dispose of the information $\U\in L^4(\Rt)$ we can derive now the identity $\U=0$ as follows: multiplying equation (\ref{N-S-f-nulle-intro})  by $\U$ and integrating on the whole space $\Rt$ we have 
$$ \int_{\Rt}(-\Delta \U)\cdot \U dx= \int_{\Rt} ((\U \cdot \vec{\nabla}) \U )\cdot \U dx + \int_{\Rt} \vec{\nabla} P  \cdot \U dx,$$ where due to the fact  $\U \in \dot{H}^{1}\cap L^4(\Rt)$ each term in this identity is well-defined.  Indeed, for the term in the left-hand side remark that as $\U \in \dot{H}^{1}(\Rt)$ then  we have $-\Delta \U \in \dot{H}^{-1}(\Rt)$.  Then, for the first term in the right-hand side, as $div(\U)=0$ we write $(\U \cdot \vec{\nabla}) \U=div(\U \otimes \U)$ where, as $\U \in L^{4}(\Rt)$ by the H\"older inequalities we have $\U \otimes \U \in L^{2}(\Rt)$ and then $div(\U \otimes \U)\in \dot{H}^{-1}(\Rt)$. Finally, in order to study the second term in the right-hand side we write the pressure $P$ as $P= \frac{1}{-\Delta} div (div (\U \otimes \U))$ hence we get $P \in L^{2}(\Rt)$ (since we have $\U \otimes \U \in L^{2}(\Rt)$) and then $\vec{\nabla} P \in \dot{H}^{-1}(\Rt)$. \\
\\
Now, integrating by parts each term in the identity above we have that $\int_{\Rt}(-\Delta \U)\cdot \U dx = \int_{\Rt}\vert \vec{\nabla} \otimes \U \vert^2 dx$, and moreover $\int_{\Rt} ((\U \cdot \vec{\nabla}) \U )\cdot \U dx=0$ and $\int_{\Rt} \vec{\nabla} P  \cdot \U dx=0$. With this identities we get $\int_{\Rt}\vert \vec{\nabla} \otimes \U \vert^2 dx =0$ and thus we have $\U=0$. 	\finpv 



\end{document}